\documentstyle[epsfig,11pt]{article}
\onecolumn 
\begin{document}

%\input{symbols}
%for the meaning of the symbols see notation.tex
\newcommand {\Pind} {\hspace{.51cm}}  % paragraph indent 
\newcommand {\Dminr}{\bigg(\frac{Vol(b^{-1}(r_1)}{2^n}\bigg)^{(1/(n-1))}}
\newcommand {\Dmin}{\bigg(\frac{V_\infty/2}{2^n}\bigg)^{(1/(n-1))}}
\newcommand {\Om} {\Omega}
\newcommand {\vare}{\varepsilon}
\newcommand {\grad}{\bigtriangledown}
\newcommand {\Ftil}{\tilde{{\cal{F}}}}  %Get from ChCo
\newcommand {\F}{{\cal{F}}}%Get from ChCo
\newcommand {\Gtil}{\tilde{{\cal{G}}}}  %Get from ChCo
\newcommand {\G}{{\cal{G}}}%Get from ChCo
\newcommand {\RR}{\,\hbox{\rm R}\!\!\!\!\!\hbox{\rm I}\,\,\,}
\newcommand {\NN}{\,\hbox{\rm N}\!\!\!\!\!\hbox{\rm I}\,\,\,}
\newcommand {\implies}{Longrightarrow}
\newcommand {\ProofEnd} {\hfill \nobreak \square \medbreak}
\newcommand{\square}{\hbox{${\vcenter{\hrule height.4pt
  \hbox{\vrule width.4pt height6pt \kern6pt
     \vrule width.4pt}
  \hrule height.4pt}}$}}

\newtheorem {theorem} {Theorem}
\newtheorem {proposition} [theorem] {Proposition}
\newtheorem {lemma} [theorem] {Lemma}
\newtheorem {definition} [theorem] {Definition}
\newtheorem {corollary} [theorem] {Corollary}
\newtheorem {remark} [theorem] {Remark}
\newtheorem {note} [theorem] {Note}
\newtheorem {example} [theorem] {Example} 

%CHANGE TITLE FORMAT!!!!

\begin{titlepage} 
\begin{center}
{\huge\bf On Loops Representing Elements\\
\vspace*{2ex} of the Fundamental Group of\\
\vspace*{2ex}  a Complete Manifold with\\ 
\vspace*{2ex} Nonnegative Ricci Curvature}\\ 
\vspace*{2ex} 

\vspace*{7ex} 
{\Large Christina Sormani} \\
\vspace*{7ex} 
{\large Johns Hopkins University} \\
\vspace*{7ex}         
{\large February, 1999} \\

\end{center}

\end{titlepage} 
\pagestyle{plain}  %try
\pagenumbering{arabic}

%\setcounter{figure}{0}
%\renewcommand{\thefigure}{\arabic{chapter}.\arabic{figure}}

%The Paper:

Throughout this paper we are concerned with the topological properties of
a complete noncompact manifold, $M^m$.  In 1969, Gromoll and Meyer proved 
that if $M^m$ has positive sectional curvature then it is diffeomeorphic 
to $\RR^m$ [GrMy].  Cheeger and Gromoll proved 
that if $M^m$ has nonnegative sectional curvature then it is diffeomorphic 
to the normal bundle over a compact totally geodesic
submanifold called the {\em soul} (1972) [ChGl2].  These theorems do not 
hold, however, if $M^m$ is only required to have nonnegative Ricci curvature.
This is demonstrated by examples of 
Nabonnand, Wei, and Sha and Yang [Nab],[Wei], [ShaYng].  Nevertheless, there
are some serious topological restrictions on such manifolds.

In 1969, Milnor conjectured that a manifold with
nonnegative Ricci curvature has a finitely generated fundamental group
[Mil].  Thus far there has been significant work done in this area
but as yet there is no counter example or proof.  Abresch-Gromoll,
Anderson, Li, and the author all have partial proofs given additional
conditions on volume or diameter  [AbrGrl] [Li] [And] [So].  This paper
concerns the properties of loops representing elements of the fundamental
group for a manifold on which there are no conditions other than those on 
Ricci curvature.   The author hopes it will prove useful to those who are 
working on this conjecture. 

Roughly, we say that $M^m$ has the {\em loops to infinity property} if 
given any noncontractible closed curve, $C$, and any compact set $K$, 
there exists a closed curve contained in $M \setminus K$ which 
is homotopic to $C$ [see Defns~\ref{defnoneloop} and~\ref{defnallloop}].  
That is, the loops can ``slide'' out to infinity beyond any compact set.  
All the examples mentioned above can be shown to have this property.  

In this paper, we prove that if $M^m$ has
positive Ricci curvature then it has the loops to infinity property
[Thm~\ref{PosLoopInf}].  As a consequence, it is impossible 
to take a manifold, remove a ball and edit in a region whose fundamental
group has more generators than the fundamental group of the boundary of the
removed ball [Thm~\ref{Boundary}].  
Thus one cannot hope to use a sequence of local surgeries to create a 
manifold with positive Ricci curvature and an infinitely generated 
fundamental group.

Furthermore, we prove that if a manifold has 
nonnegative Ricci curvature and does not satisfy the loops to infinity
property then $M^m$ is isometric to a flat normal bundle over a compact 
totally geodesic submanifold [Thm~\ref{sortasoul}].  In fact,  $M^m$ has 
a double cover which splits isometrically [Thm~\ref{doublecover}]. 

In the first section of this paper, definitions and theorems are stated
precisely with a discussion of examples and consequences.  In the
second section, we prove the Line Theorem, Theorem~\ref{noloopline}, 
which relates the loops to infinity property to the existence of a line 
in the universal cover and allows us to apply the Cheeger-Gromoll Splitting 
Theorem [ChGl1].
In the third section, the Line Theorem and a careful analysis of the
deck tranforms, leads to proofs of Theorems~\ref{doublecover}
and~\ref{sortasoul}.

Local topological consequences of the loops to 
infinity property and the previous theorems appear in the fourth section.  
Theorems~\ref{Boundary} and Theorem~\ref{doublelocal}, which are partial
extensions of theorems of Frankel [Fra], Lawson [Law], and 
Schoen-Yau [SchYau2], are included in this section.  
Further applications of these results will appear in work by the author
and Zhongmin Shen.  %ADD REF WHEN READY

The author would like to thank Professors Jack Morava and Dan Christensen
for helpful discussions regarding the topological consequences of the
loops to infinity property, and Professor Anderson for refering her to
the work of Frankel and Lawson.  Necessary background material can be found 
in texts by do Carmo, Li, and Munkrees [doC] [Li] [Mnk].

\begin{section} {Statements}

Let $M^m$ be a complete noncompact manifold.  A ray, 
$\gamma:[0,\infty)\to M$, is a minimal geodesic parametrized 
by arclength.  In contrast, a line, $\gamma:(-\infty,\infty) \to M$
is a minimal geodesic parametrized by arclength in both directions.
That is $d(\gamma(t), \gamma(s))=|t-s|$ for all $s,t \in \RR$.

A loop is a closed curve starting and ending
at a base point.  A geodesic loop is a smooth geodesic
except at the base point.

\begin{definition} {\em
Given a ray $\gamma$ and a loop $C:[0,L] \to M$ based at $\gamma(0)$,
we say that a loop $\bar{C}:[0,L] \to M$ is 
{\em homotopic to} $C$ along $\gamma$
if there exists $r>0$ with $\bar{C}(0)=\bar{C}(L)=\gamma(r)$ and the loop,  
constructed by joining $\gamma$ from $0$ to $r$ with $C$ from $0$ 
to $L$ and then with $\gamma$ from $r$ to $0$, is homotopic to $C$. }
\end{definition}

\begin{definition} \label{defnoneloop} {\em
An element $g\in \pi_1(M, \gamma(0))$, has the {\em (geodesic) loops
to infinity property along $\gamma$} if given any compact set
$K \subset M$ there exists a (geodesic) loop, $\bar{C}$,
contained in $M\backslash K$ which is
homotopic along $\gamma$ to a representative loop, $C$,
such that $g=[C]$. }
\end{definition}

It is easy to see that if $M$ is an isometrically split manifold, 
$M=N\times \RR$, with
$\gamma$ in a split direction then every $g\in \pi_1(M, \gamma(0))$
has the geodesic loops to infinity property along $\gamma$.  

\begin{example} \label{1example} {\em
If $N^3= \RR^3/G $ where $G$ is the group generated by 
$g(x,y,z)=(-y, x, z+1)$.  Then $g$ has the geodesic loops to
infinity property along $(t,0,0)$, because line segments in $\RR^3$ joining
$(t,0,0)$ to $(0,t,1)$ project to geodesic loops.  On the other hand
$g^2$ does not have the geodesic loops to
infinity property along $(t,0,0)$, because line segments in $\RR^3$ joining
$(t,0,0)$ to $(-t,0,2)$ all pass through $(0,0,1)$.  However, it does
have the loops to infinity property, because line segments from
$(t,0,0)$ to $(0,t,1)$ joined with line segments from
$(0,t,1)$ to $(-t,0,2)$ project to the required loops. }
\end{example}

\begin{example} \label{Moebius} {\em
The complete flat Moebius Strip is $M^2=\RR^2/G$ where G is generated 
by $g(x,y)=(-x, y+1)$.  Note that $g$
doesn't even have the loops to infinity property along $(t,0)=(-t,1)$
because curves joining $(t,0)$ to $(-t,1)$ must pass through 
the compact set $(\{0\}\times \RR)/G$.}
\end{example}

\begin{definition} \label{defnallloop} {\em
$M^n$ has the {\em (geodesic) loops to infinity} property if given any 
given any ray, $\gamma$, and any element $g\in \pi_1(M,\gamma(0))$,
$g$ has the (geodesic) loops to infinity property along $\gamma$.}
\end{definition}

\begin{example} \label{NabWei}  {\em
Nabannond has constructed an example
of a manifold with positive Ricci curvature 
which is diffeomorphic to $\RR^3 \times S^1$ [Nab].  Wei
has constructed examples which are diffeomorphic to $\RR^k \times N$
where the fundamental group, $\pi_1(N)$, is any torsion free
nilpotent group [Wei].  One can show that the examples of Nabonnand and 
Wei actually satisfy the {\em geodesic} loops to infinity property.}
\end{example}

In Section 2, we will prove the following theorem.

\begin{theorem} [Line Theorem] \label{noloopline} \label{line}
If $M^n$ is a complete noncompact manifold which does not
satisfy the geodesic loops to infinity property then there
is a line in its universal cover.
\end{theorem}

Recall that the Splitting Theorem of Cheeger and Gromoll
states that a manifold with nonnegative Ricci curvature
which contains a line splits isometrically [ChGl2] [see also
Li's text, Thm 4.2].  Thus we have the following consequence
of Theorem~\ref{noloopline}.

\begin{theorem} \label{PosLoopInf}
If $M^n$ is complete noncompact with $Ricci \ge 0$ and there exists
$y\in M^n$ such that $Ricci_y >0$, then $M^n$ has the geodesic loops to 
infinity property.
\end{theorem}

In Section 3, we prove more in the case of nonnegative Ricci curvature.

\begin{proposition} \label{gbad}  {\em
If $M^n$ is a complete noncompact manifold with $Ricci \ge 0$
and there exists an element $g \in \pi_1(M)$
which does not satisfy the {\em geodesic loops to infinity property}
along a given ray $\gamma$, then the universal cover splits isometrically,
$\tilde{M}= N^{m-k} \times \RR^k$.  Furthermore the lift $\tilde{\gamma}$
of $\gamma$, is in the split direction,
\begin{equation}
\tilde{\gamma}(t)=(x(0), y(t))
\end{equation}
and }
\begin{equation}
g_*(\tilde{\gamma}'(t))=-\tilde{\gamma}'(t).
\end{equation}
\end{proposition}

\begin{corollary} \label{gsquared}
If $M^n$ is a complete noncompact manifold with $Ricci \ge 0$
and $g \in \pi_1(M)$, then either $g$ or $g^2$ has the geodesic
loops to infinity property.
\end{corollary}

If we consider manifolds which don't even satisfy the weaker loops to 
infinity property, we get a stronger result.

\begin{theorem} \label{doublecover} [Double Cover Theorem]
If $M^n$ is a complete noncompact manifold with $Ricci \ge 0$
and there exists an element $g \in \pi_1(M)$
which does not satisfy the loops to infinity property
along a given ray $\gamma$, then all elements
$h \in \pi_1(M, \gamma(0))$ satisfy 
\begin{equation} \label{allh}
h_*(\tilde{\gamma}'(t))= \pm \tilde{\gamma}'(t).
\end{equation}
Furthermore, $M^m$ has a split double cover which lifts $\gamma$ to a line.
\end{theorem}

Local consequences of this theorem appear in Section 4.
[Theorem~\ref{doublelocal}].

Cheeger and Gromoll proved that any complete noncompact manifold with
nonnegative sectional curvature is diffeomorphic to a normal bundle
over a compact totally geodesic submanifold called a soul [ChGr2].  This
is not true for manifolds with nonnegative Ricci curvature [ShaYng].

However, $M$ does  have a soul if it doesn't have the geodesic loops 
to infinity
property along any ray.  This soul is defined using Busemann functions,
which are reviewed in Section~\ref{LineThm} above Lemma~\ref{buseback}.

\begin{theorem} \label{sortasoul}
If $M^n$ is a complete noncompact manifold with $Ricci \ge 0$
and there exists an element $g \in \pi_1(M)$
which does not satisfy the loops to infinity property
along a given ray $\gamma$, then the Busemann function, $b_\gamma$ associated
with that ray has a minimum
$$
-s_\gamma=min_{x\in M} (b_\gamma(x))
$$ 
and $M^n$ is a flat normal bundle over $b_\gamma^{-1}(-s_\gamma)$.

If $g\in \pi_1(M)$ doesn't satisfy the geodesic loops to infinity property
along any ray, $\gamma$, then $M^n$ is a flat normal bundle
over a compact totallt geodesic soul, $S$, where
$S=\bigcap_\gamma(b_\gamma^{-1}(-s_\gamma)).$
\end{theorem}

This is the same soul as the one produced in Cheeger and Gromoll's paper
if $M^m$ has nonnegative sectional curvature.

Note that in Example~\ref{1example}, the soul, $S=\{(0,0,z): z\in [0,1]\}$.

In general, for manifolds with nonnegative Ricci curvature, it is an 
open question whether, $\bigcap_{\gamma}(b_{\gamma}^{-1}(-s_\gamma))$
is compact or not.   

In Section 4, we discuss some topological consequences of the loops
to infinity property.  In particular Theorem~\ref{Boundary}, states that in
a manifold with the loops to infinity property, the group homomorphism
induced by the inclusion from $\pi_1(\partial D) \to  \pi(Cl(D))$ is a 
surjection.
Thus, if the boundary of a region in a manifold with positive Ricci
curvature is simply connected, the region must
be simply connected as well.  This consequence is an old theorem of Schoen
and Yau [SchYau2].   A similar weaker theorem is proven if
$M^m$ has $Ricci \ge 0$ [Theorem~\ref{doublelocal}].

\end{section}

\begin{section} {Proof of the Line Theorem} \label{LineThm}

\Pind
In this section, $M$ is a Riemannian Manifold and we make no assumptions
on Ricci curvature.
We begin with a construction of the line in the universal cover.
Elements of this proof are used again to prove other theorems in 
the next section.

{\bf Proof of Theorem~\ref{noloopline}:}
Let $\gamma$ be a ray, $g\in \pi_1(M, \gamma(0))$ such that $g$ doesn't
satisfy the geodesic loops to infinity property.  Let $C$ be a 
representative of $g$ based at $\gamma(0)$.  There exists a compact set $K$, 
such that there is no closed geodesic
contained in $M\backslash K$ which is
homotopic to $C$ along $\gamma$.  Let $R_0 >0$ such that
$B_{x_0}(R_0)\supset K$.

So for all $r >R_0$,
any loop based at $\gamma(r)$ which is homotopic to $C$ along $\gamma$
must pass through $K$.  Let $r_i>R_0$ be an increasing sequence diverging
to infinity.

Let $\tilde M$ be the universal cover and $\tilde{C}$ be a lift
of $C$ running from $\tilde{x_0}$ to $g\tilde{x_0}$.  Since $C$ is noncontractible,
$g$ is not the identity and $\tilde{x_0}\neq g\tilde{x_0}$.  
Let $\tilde{\gamma}$ be the lift of $\gamma$ starting at $\tilde{x_0}$
and $g\tilde{\gamma}$ be the lift  starting at $g\tilde{x_0}$.
Then if $\tilde{C_i}$ is a minimal geodesic from
$\tilde{\gamma}{(r_i)}$ to $g\tilde{\gamma}{(r_i)}$, $C_i=\pi \circ \tilde{C_i}$
is a loop based at $\gamma(r_i)$ which is homotopic to $C$ along $\gamma$.
Thus there exists $t_i$ such that $C(t_i)\subset K$.

Let $L_i=L(C)=d_{\tilde{M}}(\tilde{\gamma}{(r_i)},g\tilde{\gamma}{(r_i)})$.  

Let $\tilde{K}$ be the lift of $K$ to the fundamental domain in
$\tilde{M}$ such that $\tilde{x_0}\in\tilde{K}$. Note that
$\tilde{K}$ is precompact.

 For all $i\in \NN$ there is an element $g_i \in \pi_1(M, x_0)$ such that
$g_i \tilde{C}(t_i) \in \tilde{K}$.

Note that 
$$
t_i= d_{\tilde{M}}(\tilde{\gamma}{(r_i)},\tilde{C}(t_i)) \ge 
 d_{{M}}({\gamma(r_i)},{C}(t_i))
\ge d_M(\gamma(r_i), K) \ge r_i-R_0
$$
and
$$
L_i-t_i=d_{\tilde{M}}(g\tilde{\gamma}{(r_i)},\tilde{C}(t_i)) \ge 
 d_{{M}}({\gamma(r_i)},{C}(t_i))
\ge d_M(\gamma(r_i), K) \ge r_i-R_0.
$$

So $g_i \tilde{C}$ are minimal geodesics running from
$(t_i-(r_i-R_0))$ to $(t_i+(r_i-R_0))$ such that 
$g_i \tilde{C}(t_i) \in \tilde{K}$.  Taking $r_i$ to infinity,
a subsequence of $( g_{i\,*} \tilde{C}'(t_i))$ must 
converge to a unit vector $(\gamma_\infty'(0))$ based 
at $\gamma_\infty(0)\in Cl(\tilde{K})$. Furthermore, the geodesic, 
\begin{equation}\label{herealine}
\gamma_\infty(t)=\exp_{\gamma_\infty(0)}(t\gamma_\infty'(0))
\end{equation}
is a line.

\ProofEnd

\begin{note} \label{Lirinote}
Note that $\lim_{i\to\infty} L_i/r_i =2$ because
$L_i \le 2r_i + d(\tilde{\gamma}(0),g\tilde{\gamma}(0))$
and $L_i \ge L_i-t_i +t_i\ge 2(r_i-R_0)$.
\end{note} 

Recall that given a ray, $\gamma$, parametrized by arclength,
then the Busemann function associated with that ray, 
$b_\gamma: M \to \infty$ is defined,
$$
b_\gamma(x)=\lim_{R\to\infty} R-d_{M}(x, \gamma(R)).
$$
For example, if $M$ is Euclidean space, then 
$\gamma(t)=\gamma(0)+ t \gamma'(0)$,
and
\begin{equation} \label{beuchlid}
b_\gamma(x)=<x-\gamma(0), \gamma'(0)>.
\end{equation}
If $M$ is a manifold with nonnegative Ricci curvature that contains a line,
$\gamma$, then by Cheeger and Gromoll, $M$ is the isometric product
of $\gamma(\RR)$ and $b_\gamma^{-1}(\{0\})$.

The following lemma is useful in analyzing the properties of deck
transforms in conjunction with rays.  It will be used in the next section.

\begin{lemma} \label{buseback}
If $\tilde{\gamma}$ is the lift of a ray $\gamma$ then 
for all deck transforms $g$,
$$
b_{\tilde{\gamma}}(g\tilde{\gamma}(a)) \le a.
$$
\end{lemma}

\noindent {\bf Proof:}
For any $x \in \tilde{M}$ we have
$$
d_{\tilde{M}}(gx, \tilde\gamma(R)) \ge d_M(\pi(x), \gamma(R)).
$$
If we subtract $R$ on both sides and take R to infinity,
we get
$$
-b_{\tilde{\gamma}} (gx)\ge -b_{\gamma} (\pi(x)).
$$
Setting $x=\tilde{\gamma}(a)$, then $\pi(x)=\gamma(a)$ and we are done.
\ProofEnd

\end{section}

\begin{section} {Nonnegative Ricci Curvature}

Throughout this section we assume that $M^m$ is a manifold
with nonnegative Ricci curvature which does not satisfy the geodesic loops
to infinity property.  Thus by the Line Theorem, its universal
cover, $\tilde{M}$ contains a line.  So by the Splitting Theorem of Cheeger 
and Gromoll, the universal cover splits isometrically 
into $N^{m-k}\times \RR^k$ where $N^{n-k}$ has no lines and $k \ge 1$.

Let $p_{\RR}:\tilde{M} \to \RR^k$ and
$p_N:\tilde{M}\to N$ be the projections.  If $g:\tilde{M}\to\tilde{M}$ 
is an isometry, then it acts on each component seperately. [REF Cheeger] 
That is $g=g|_N:N\to N$ and $g=g|_{\RR^k}:\RR^k \to \RR^k$ are isometries.  

The following lemma is quite easy to prove.

\begin{lemma} \label{projmin}
  If $\eta:(a,b)\to \tilde{M}$ is minimal then 
$p_{\RR} \circ \eta$ and $p_N \circ \eta$ are minimal geodesics
between their endpoints.  It is possible that one of them is constant.
\end{lemma}

Note, however, that the geodesics in this lemma are paremetrized proportional
to arclength and are not normalized like the geodesics, rays and lines
constructed in the proof of Theorem~\ref{noloopline}. 

We first prove Proposition~\ref{gbad}.  There is a given ray
$\gamma$ and a given element
$g \in \pi_1(M, \gamma(0))$ which does not satisfy the 
{\em geodesic loops to infinity 
property} along  $\gamma$. We must show that $\gamma$ lifts to
the purely split direction and that 
$g_*(\tilde{\gamma}'(t))=-\tilde{\gamma}'(t)$.

\noindent{\bf Proof of Proposition~\ref{gbad}:}
By the proof of the Line Theorem, we know there are minimal geodesics
$\tilde{C_i}$, running from $\tilde{\gamma}(r_i)$
to $g\tilde{\gamma}(r_i)$, whose lengths $L_i$, are growing like $2r_i$.
So, intuitively,  $\tilde{\gamma}$ and $g\tilde{\gamma}$ should be in
the opposite directions, and thus can only fit in the split direction.

Let $\tilde{C_i}$ be the curves constructed in Theorem~\ref{noloopline}.
Let $x_i(t)=p_N(\tilde{C_i}(t))$ and $y_i(t)=p_{\RR}(\tilde{C_i}(t))$.
By Lemma~\ref{projmin}, these are minimal geodesics from $[0, L_i]$.
Since $y_i$ is a minimal geodesic in Euclidean space, it can
be written as 
$$
y_i(t)= y_i'(t_i)(t-t_i)+y_i(t_i)
$$
where $t_i \in (0,L_i)$ as in Theorem~\ref{noloopline}.  Since all minimal 
geods are parametrized proportional to arclength, $|x_i'(t)|=|x_i'(t_i)|$.
Since  $\tilde{C_i}$ are minimal geodesics parametrized by arclength,
$$
|x_i'(t_i)|^2 +|y_i'(t_i)|^2 =|\tilde{C_i}'(t_i)|^2=1.
$$ 

Now, $g_{i\,*}(\tilde{C_i}'(t_i))$ converges to 
$\gamma_\infty'(0)$, where $\gamma_\infty$ is a line defined in 
(\ref{herealine}).  Let
$x_\infty(t)=p_N(\gamma_\infty(t))$ and 
$y_\infty(t)=p_{\RR}(\gamma_\infty(t))$.
By Lemma~\ref{projmin}, $x_\infty(t)$ and $y_\infty(t)$ are lines or
constants.  Since $N$ contains no lines, $x_\infty(t)$ is a constant.
Thus 
$$
|y_\infty'(0)|=|\gamma_\infty'(0)|=1.
$$  

Since $g_i|_N$ is an
isometry,
$$
\lim_{i\to\infty} |x_i'(t_i)| =
\lim_{i\to\infty}|g_{i\,*} x_i'(t_i)| =|x_\infty'(0)| =0.
$$
Similarly,
\begin{equation} \label{yinfty} 
\lim_{i\to\infty} |y_i'(t_i)| =\lim_{i\to\infty}|g_{i\,*} y_i'(t_i)|
=|y_\infty'(0)|= 1. 
\end{equation}

Let $\tilde{\gamma}$ be the lifted ray as in Theorem~\ref{noloopline}.
Let $x=p_N\circ \tilde{\gamma}$ and $y=p_{\RR}\circ \tilde{\gamma}$.
Recall that $\tilde{C_i}$ is a minimal geodesic of length $L_i$ from
$\tilde{\gamma}(r_i)$ to $g\tilde{\gamma}(r_i)$.  
So $x_i$ is minimal from $x(r_i)$ to $g(x(r_i))$ and
$y_i$ is minimal from $y(r_i)$ to $ g(y(r_i))$.
Thus
$$
d_{\RR^k}\big(y(r_i), g(y(r_i))\big)= L_i |y_i'(t_i)|.
$$
Taking $i\to \infty$, and applying (\ref{yinfty}), we have
\begin{equation} \label{2yinfty}
\lim_{i\to\infty}\frac{d_{\RR^k}(y(r_i), g(y(r_i))}{L_i}=1.
\end{equation}

Now $y(t)$ is a minimal geodesic in $\RR^k$, so $y(t)=y'(0)t+y(0)$
and $g(y(t))$ is also a minimal geodesic, so $g(y(t))=g_*y'(0)t+g(y(0))$.
Thus
\begin{eqnarray}
\frac{d_{\RR^k}(y(r_i), g(y(r_i))}{L_i}&=&
\frac{ |y'(0)r_i+y(0)-(g_*y'(0)r_i+g(y(0)))|}{L_i}\\
&\le & \frac{ |y'(0)-g_*y'(0)|r_i +|y(0)- g(y(0)))|}{L_i}.
\end{eqnarray}
Putting this together with (\ref{2yinfty}), we have 
\begin{eqnarray}
\lim_{i\to \infty}\frac{ |y'(0)-g_*y'(0)|r_i +|y(0)- g(y(0)))|}{L_i}=1.
\end{eqnarray}
So $|y'(0)-g_*y'(0)|\neq 0$.  In particular, $y'(0) \neq 0$, and the
original lifted ray $\tilde{\gamma}$ has a component in the split direction.

By Note~\ref{Lirinote}, 
$$
\lim r_i/L_i =1/2,
$$ 
so
$$
|y'(0)-g_*y'(0)|=2.
$$
However 
$$
|y'(0)| \le |\tilde{\gamma}'(0)|=1.
$$
Thus
\begin{equation} \label{1doubleg}
2 =|y'(0)-g_*y'(0)|\le |y'(0)|+| g_*y'(0)|=2|y'(0)|\le 2
\end{equation}
and $y'(0)=\tilde{\gamma}'(0)$.  So the original lifted ray is completely
in the split direction.  

Furthermore, to have equalities in (\ref{1doubleg}),
the element $g \in \pi_1(M)$ which did 
not have geodesic loops to infinity must satisy 
\begin{equation} \label{doubleg}
g_*\tilde{\gamma}'(0)=g_*y'(0)=-y'(0)=-\tilde{\gamma}'(0).
\end{equation}

\ProofEnd

We've completed the proof of Proposition~\ref{gbad} and 
Corollary~\ref{gsquared} follows trivially.

We now turn to a manifold with nonnegative Ricci curvature which does
not even have the loops to infinity property.  That is, we have a ray $\gamma$
and an element $g\in \pi_1(M, \gamma(0))$ which does not have the loops to
infinity property along $\gamma$.  Thus, $g$ also does not have the
geodesic loops to infinity property.  

\noindent {\bf Proof of Theorem~\ref{doublecover}:}
All the conditions of Theorem~\ref{gbad} hold so we know that $\tilde{\gamma}$
lifts to a the split direction of the universal cover and (\ref{doubleg})
holds.

We first try to construct loops to infinity which may not be geodesic loops.  
Keep in mind that if $M$ is the moebius strip then there are no such loops
while in Example~\ref{1example}, there are.

Now we've shown that $\tilde{\gamma}:[0,\infty)\to \tilde{M}$ lies 
completely in the split direction, so in fact it extends to a line
$\tilde{\gamma}:(-\infty, \infty)\to \tilde{M}$.  If this line
projects to a line $\gamma=\pi\circ \tilde{\gamma}:(-\infty, \infty)\to M$ 
then $M$ splits and again we have geodesic loops to infinity.  So
it does not project to a line.

Suppose there exists $h \in \pi_1(M)$ such that 
\begin{equation}\label{assumph}
h_*(\tilde{\gamma}'(0))\neq -\tilde{\gamma}'(0)
\textrm{ and } h_*(\tilde{\gamma}'(0))\neq \tilde{\gamma}'(0).
\end{equation}

Then $h$ has geodesic loops to infinity.  That is, for all $R>0$
there exist $r_i >R$ such that the minimal geodesic, $\eta_i$, of
length, $l_i$, from $\tilde{\gamma}(r_i)$ to $h\tilde{\gamma}(r_i)$, satisfy
$$
\pi(\eta_i([0,l_i])) \cap B_{\gamma(0)}(R) = \emptyset.
$$

Furthermore, by the fact that $h_* \in O(n)$ and ((\ref{assumph}) holds, 
$$
(h^{-1}g)_*(\tilde{\gamma}'(0))=h^{-1}_*(-\tilde{\gamma}'(0))
= -h^{-1}_*(\tilde{\gamma}'(0)) /ne -\tilde{\gamma}'(0).
$$
Thus $h^{-1}g$ has the geodesic loops to infinity property, and 
for all $R>0$
there exist $\bar{r}_i >R$ such that the minimal geodesic, $\bar{\eta}_i$, of
length, $\bar{l}_i$, from $\tilde{\gamma}(\bar{r}_i)$ to 
$h^{-1}g\tilde{\gamma}(\bar{r}_i)$, satisfy
$$
\pi(\bar{\eta}_i([0,\bar{l}_i])) \cap B_{\gamma(0)}(R) = \emptyset.
$$ 
Note that $h\bar{\eta}_i$ is a minimal geodesic from
$h\tilde{\gamma}(\bar{r}_i)$ to $g\tilde{\gamma}(\bar{r}_i)$
such that
$$
\pi(h\bar{\eta}_i([0,\bar{l}_i])) \cap B_{\gamma(0)}(R) = \emptyset.
$$

Thus for all $R>0$ there exists an $r_i>0$, such that there is a curve, 
$s:[0, l_i +\bar{l_i}+ 2|r_i-\bar{r}_i|] \to \tilde{M}$
running from $\tilde{\gamma}(r_i)$ to $g\tilde{\gamma}(r_i)$.  This curve, $s$,
first runs along $\eta_i$ from $\tilde{\gamma}({r}_i)$ to 
$h{\gamma}({r}_i)$, then along 
$h\tilde{\gamma}(t)$ to $h\tilde{\gamma}(\bar{r}_i)$,
then along $h\bar{\eta_i}$ to $g\tilde{\gamma}(\bar{r}_i)$, and
finally along $g\tilde{\gamma}(t)$ to $g\tilde{\gamma}(r_i)$.
Clearly $s(t)=\pi(\tilde{s}(t))$ avoids $B_{\gamma(0)}(R)$
and is homotopic along $\gamma$ to any curve based at $\gamma(0)$
representing $g$.

This contradicts the hypthesis in (\ref{assumph}), so we have proven
(~\ref{allh}). 

We now construct the double cover.  Let 
$$
H =\{h \in \pi_1(M, \gamma(0)): 
h_*(\tilde{\gamma}'(0))=\tilde{\gamma}'(0)\}.
$$
Then $H$ is clearly a normal subgroup of $\pi_1(M, \gamma(0))$
and $\pi_1(M, \gamma(0))/H=\{[e],[g]\}$.
Thus there exists a double cover $\bar{M} =\tilde{M}/H$ of $M$.

Let $\pi_H:\tilde{M} \to \bar{M}$.  We claim that $\pi_H(\tilde{\gamma})$
is a line.

Suppose not.  Then there exists $s>0$ such that 
$$
d_{\bar{M}}(\pi_H(\tilde{\gamma})(-s),\pi_H(\tilde{\gamma})(s)) <2s.
$$
So there exists $\bar{h}\in H$, such that
\begin{equation} \label{notline2}
d_{\tilde{M}}(\bar{h}(\tilde{\gamma}(-s)),\tilde{\gamma}(s)) <2s.
\end{equation}

Let $(\bar{x_0}, \bar{y_0})=\bar{h}(\tilde{\gamma}(0))$. 
Now by Lemma~\ref{buseback}, 
\begin{equation} \label{bhere}
b_{\tilde{\gamma}}((\bar{x_0}, \bar{y_0})) \le b_{\tilde{\gamma}}(\gamma(0)).
\end{equation}
Since $\tilde{M}$ is split and $\tilde{\gamma}$ is in the split
direction,
$$
b_{\tilde{\gamma}}(x_0, y(0))=b_{\tilde{\gamma}}(x(0), y(0))
$$
and as in (\ref{beuchlid}),
\begin{equation} \label{1busebackap}
b_{\tilde{\gamma}}((x(0), y(0)+v))=b_{\tilde{\gamma}}((x(0), y(0))) 
+<v, y'(0)>_{\RR^k}.
\end{equation}
Setting $v=\bar{y}(0) -y(0)$ and applying (\ref{bhere}), we have
\begin{eqnarray} \label{binner}
<\bar{y}_0-y(0), y'(0)>& = & b_{\tilde{\gamma}}(x(0), \bar{y}_0)
                             - b_{\tilde{\gamma}}(x(0), y(0)) \\
&=& b_{\tilde{\gamma}}(\bar{x}_0, \bar{y}_0)
                             - b_{\tilde{\gamma}}(x(0), y(0)) \le 0 \\
\end{eqnarray}

Now $h_*(\tilde{\gamma}'(0))=\tilde{\gamma}'(0)$, so
$h(\tilde{\gamma}(t))=(\bar{x_0}, -y'(0)t + \bar{y_0})$ while \newline
$\tilde{\gamma}(t)= ({x(0)}, -y'(0)t + {y(0)})$.
Thus, by (\ref{binner}),
\begin{eqnarray*} \label{notline3}
d_{\tilde{M}}(\bar{h}(\tilde{\gamma}(-s)),\tilde{\gamma}(s))^2 
&\ge& d_{\RR^k}(-y'(0)(-s) + \bar{y_0}, -y'(0)s + {y(0)})^2\\
&=& |2s y'(0) + \bar{y_0}-y(0)|^2\\
&=& |2s y'(0)|^2 + 4s<\bar{y_0}-y(0),y'(0)> +|\bar{y_0}-y(0)|^2 \\
&\ge & |2s y'(0)|^2.
\end{eqnarray*}
This contradicts (\ref{notline2}), so $\pi_H(\tilde{\gamma})$ is a line
in the double cover.  

% LOOK OVER
%Furthermore, all the inequalities in (\ref{notline3}) are equalities,
%so $h(y(0)=\bar{y_0}=y(0)$.  In fact, for all $h\in H$, viewed as isometries
%on $\RR^k$, $h(y(0)))=y(0)$, so $h(y(t))=y(t) \forall t\in\RR$. 

\ProofEnd

We now prove Theorem~\ref{sortasoul}, in which we study the Busemann
functions on manifolds which don't satisfy the loops to infinity
property.

\noindent {\bf Proof of Theorem~\ref{sortasoul}:}

Now $\gamma(t)$ is not a line, else we would have had geodesic rays to
infinity.  However, it is possible that $\gamma:[-s, \infty) \to M$ is
a ray for some $-s<0$.  Let $s_\gamma>0$, be defined such that
$\gamma:[-s_\gamma, \infty) \to M$ is a ray and $\gamma:[-s, \infty) \to M$ is
not a ray for any $s>s_\gamma$.  

We claim that any element $h \in \pi_1(M, \gamma(0))$ maps the
level set 
$$
b_{\tilde{\gamma}}^{-1}(-s_\gamma)=N^{n-k}\times b_{y}^{-1}(-s_\gamma) \subset\tilde{M}
$$ 
to itself.  Recall that $\tilde{\gamma}(t)=(x(0),y(t))$ by 
Theorem~\ref{gbad}.

Recall that by (\ref{beuchlid}) and the splitting, if $(z,w)\in \tilde{M}$,
then 
\begin{eqnarray} \label{2beuchlid}
b_{\tilde{\gamma}}((z,w)) &=& b_{y}(w)=<w-y(0), y'(0)> \\ 
&=&<w-y(-s_\gamma), y'(0) > + <y(-s_\gamma)-y(0), y'(0) > \\
&=& <w-y(-s_\gamma), y'(0) > - s_\gamma.
\end{eqnarray}

Suppose $(z_0,w_0)\in b_{\tilde{\gamma}}^{-1}(-s_\gamma)$.
Since $h$ preserves the splitting, and satisfies (\ref{allh}),
\begin{eqnarray*}
b_{\tilde{\gamma}}(h(z_0,w_0)) & = &b_{\tilde{\gamma}}((hz_0,hw_0)) \\
&=& <hw_0-y(-s_\gamma), y'(0) > - s_\gamma\\
&=& <w_0-h^{-1}y(-s_\gamma), h^{-1}_*y'(0) > - s_\gamma\\
&=& \pm <w_0-h^{-1}y(-s_\gamma), y'(0) > - s_\gamma.
\end{eqnarray*}
Thus $h(z_0,w_0)\in b_{\tilde{\gamma}}^{-1}(-s_\gamma)$ iff
$$
<h^{-1}y(-s_\gamma)-y(-s_\gamma), y'(0) >=0.
$$

So we need only show that $h(\tilde{\gamma}(-s_\gamma))\in 
b_{\tilde{\gamma}}^{-1}(-s_\gamma)$.

Since $\tilde{\gamma}$ is a lift of a ray, we can apply
Lemma~\ref{buseback}.   Thus
$$
h(\tilde{\gamma}(-s_\gamma))\in b_{\tilde{\gamma}}^{-1}((-\infty,-s_\gamma]).
$$

Now $\gamma:[-s, \infty) \to M$ is not a ray for if $s>s_\gamma$.  So
there exists $s_i \to  s_\gamma$ and $r_i\to \infty$
such that
\begin{eqnarray}\label{buseless}
r_i +s_i &>& d_{M}(\gamma(-s_i), \gamma(r_i)) \\
&\ge& d_{\tilde{M}}(h \tilde{\gamma}(-s_i),\tilde{\gamma}(r_i)).
\end{eqnarray}   
Subtracting $r_i$ and taking a limit as $r_i$ approaches infinity,
\begin{eqnarray}
s_\gamma= \lim_{i\to\infty} s_i &\ge& 
\lim_{i\to\infty}(d_{M}(\gamma(-s_i), \gamma(r_i))-r_i)\\
&\ge& -\lim_{i\to\infty}(r_i-
d_{\tilde{M}}(h \tilde{\gamma}(-s_i),\tilde{\gamma}(r_i))\\
&= & -\lim_{i\to\infty}(r_i-
d_{\tilde{M}}(h \tilde{\gamma}(-s_\gamma),\tilde{\gamma}(r_i))\\
&=& -b_{\tilde{\gamma}}(h \tilde{\gamma}(-s_\gamma)).
\end{eqnarray}   
Thus
$$
h(\tilde{\gamma}(-s_\gamma)) \in b_{\tilde{\gamma}}^{-1}([-s_\gamma, \infty),
$$
and the claim is proven.

Thus $\tilde{M}$ is a flat normal bundle over $b_{\tilde{\gamma}}^{-1}(-s_\gamma)$
with one dimensional fibres, and $\pi_1(M, \gamma(0))$
is a group which preserves the base and
maps fibres to fibres.  Thus $M$ is a flat normal bundle over
$\pi(b_{\tilde{\gamma}}^{-1}(-s_\gamma))$.

Furthermore $g\tilde{\gamma}\in b_{\tilde{\gamma}}^{-1}(s_\gamma)$ and
$g\tilde{\gamma}'(0)=-\tilde{\gamma}'(0)$ implies that
$$
b_{g\tilde{\gamma}}(\tilde{p})-s_\gamma
=-(b_{\tilde{\gamma}}(\tilde{p})-s_\gamma) \qquad \forall \tilde{p}\in\tilde{M}.
$$
Thus for any $p\in M$, with lift $\tilde{p}$,
\begin{eqnarray*}
b_{{\gamma}}(p) 
&= &\lim_{R\to\infty} R-d_M(p, \gamma(R))\\
&\ge& \lim_{R\to\infty} R- 
max \{d_{\tilde{M}}(\tilde{p}, \tilde{\gamma}(R)).
      d_{\tilde{M}}(g\tilde{p}, \tilde{\gamma}(R))\} \\
&\ge & min \{b_{\tilde{\gamma}}(\tilde{p}),b_{g\tilde{\gamma}}(\tilde{p})\}\\
&\ge -s_\gamma,
\end{eqnarray*}
Thus $-s_\gamma=\min_{p\in M} b_{{\gamma}}(p)$.

Now suppose there is an element $g\in \pi_1(M,x_0)$,
such that for any ray $\gamma$ with $\gamma(0)=x_0$, $g$ doesn't have the 
geodesic loop to infinity property.  Then for each $\gamma$, we
have $s_\gamma$ and splitting in the $\gamma$ direction such that
$\tilde{M}$ is a flat normal bundle over a totally geodesic set
$b_{\tilde{\gamma}}^{-1}(-s_\gamma)$
and $\pi_1(M, \gamma(0))$
is a group which preserves the base and
maps fibres to fibres.  

Thus $\bigcap_\gamma b_{\tilde{\gamma}}{-1}(-s_\gamma))$ is totally geodesic
and 
$$
\tilde{M}=\bigcap_\gamma b_{\tilde{\gamma}}^{-1}(-s_\gamma)) \times \RR^{l}
$$
where 
$$
l=dim(span\{\gamma'(0): \gamma \textrm{ is a ray based at } x_0\}.
$$
Thus $M$ is a flat normal bundle over the totally geodesic
$S=\bigcap_\gamma b_{\gamma}^{-1}(-s_\gamma))$.

Now $S$ is totally geodesic, so if it were noncompact it would contain
a ray $\gamma$.  However, no ray is ever contained in its own level set,
so no ray can be contained in $S$.
\ProofEnd

\end{section}

\begin{section}{Topological Consequences of Loops to Infinity}

The simplest consequence of the loops to infinity property is the
following simple theorem.

\begin{theorem} \label{Pi1onto}
If $M^n$ has the loops to infinity property, $K$ is a compact set
and $y_0$ is a point in an unbounded component $U \subset M\setminus K$,
then the inclusion map 
$$
i_*:\pi_1(U, y_0) \longmapsto \pi_1(M, y_0)
$$
is onto.
\end{theorem}

\noindent {\bf Proof:}
 Since $U$ is unbounded, there exists $R_0>0$ and  a ray, $\gamma$, 
from $y_0$ such that $\gamma(r)\in U$ for all $r\ge R_0$.  

Given $g \in \pi_1(M, y_0)$, $M$ has the loops to infinity property,
so there exists a loop $\bar{C}$ contained in $M\setminus K$ which is homotopic
along $\gamma$ to a representative loop $C$ such that $[C]=g$.  
Since $U$ is a connected component, $\bar{C}\in U$.  Now we can add
segments of $\gamma$ to the front and back of $\bar{C}$ to get a curve
$\eta$ which is
homotopic to $C$, based at $y_0$ and still contained in $U$.
The $[\eta]\in \pi_1(U, y_0)$ and $i_*([\eta])=g$.
\ProofEnd

The following theorem is a localization of the above and is
proven below the statements of its corollaries.

\begin {theorem} \label{Boundary}  %put in intro
Let $M^n$ be a complete riemannian manifold with the loops to infinity
property along some ray, $\gamma$.  Let $D\subset M$ be a precompact region
with smooth boundary
containing $\gamma(t_0)$ and $S$ be any  connected component of $\partial D$
containing a point $\gamma(t_1)$.  Then the inclusion map 
$$
i_*:\pi_1(S, \gamma(t_1)) \longmapsto \pi_1(Cl(D), \gamma(t_1))
$$
is onto.
\end{theorem}

Note that this theorem is closely related to results of Frankel, 
Lawson and Schoen-Yau which concern precompact regions in a complete
manifold, $M^m$, with nonnegative Ricci curvature but without the assumption 
that $M^m$ is noncompact.   Frankel and Lawson are able to prove that
$i*$ is surjective, but they require that the boundary have conditions
on its mean curvature [Fra] [Law].  Schoen and Yau do not require any
extra boundary conditions but they have a weaker conclusion than the
one in Theorem~\ref{Boundary} [SchYau2].
%(p345 Lectures on Harm Maps)  
The methods used to prove the above theorems involve Synge's second
variation of arclength in  [Fra] and [Law], and harmonic maps in
 [SchYau2].

\begin{corollary}
Let $M^n$ has the loops to infinity property.
If $\partial D$ is simply connected, then $\pi_1(D)$ is trivial.
\end{corollary}

\begin{corollary}
If $M^n$ has nonnegative Ricci curvature and $\partial D$ is simply connected,
then $\pi_1(D)$ can only contain elements of order 2.
\end{corollary}

%Now the fundamental group of a compact manifold with nonpositive 
%curvature does not contain any elements of order 2 (WHY?????? CHECKIFTRUE).
%Thus there is no
%nontrivial homomorphism from $\pi_1(D)$ to the fundamental group 
%of a compact manifold with nonpositive curvature.  This is the theorem of
%Schoen and Yau (p345 Lectures on Harm Maps).

\begin{corollary}
Any Riemannian manifold $M^n$ with the loops to infinity property
which is simply connected at infinity, is simply connected. 
\end{corollary}

\noindent {\bf{Proof of Theorem~\ref{Boundary}:}}
Since $S$ is smooth and compact, there exists 
$$
r_0=\min_{x\in\partial D} \,\,injrad(x)\,\, >\,\,0,
$$ 
such that
the tubular neighborhood $T_{r_0}(S)$ is homotopic to $S$.
%TRUE FOR A TOP MANIFOLD?
The exponential map along the normals can be used to create the homotopy. 

Let $U=D \cup T_{r_0}(S)$ and $V= (M \,\,\setminus D)\cup T_{r_0}(S)$.
Then $U\cap V= T_{r_0}(S)$.  Note that $U$ is homotopic to $D$.  We wish
to prove that:
$$
i: \pi(T_{r_0}(S), \gamma(t_0)) \longmapsto \pi(U, \gamma(t_0))
$$
is onto.  That is, we must show that given any loop, $C_1\in U$, based at 
$\gamma(t_0) $ which is not contractible in $U$, there exists a curve $C_2
\in U\cap V$ based at the same point, which is homotopic to $C_1$.

Let 
\begin{equation}\label{t2}
t_2=\sup\{t\,\, s. t.\,\, \gamma(t)\in D\}
\end{equation}

Fix $C_1 \in U$ as above.  If $C_1$ is not contractible in $M$ then
by the loops to infinity property and the 
compactness of $Cl(U)$, there exists a loop $C_3 \in M\,\,\setminus Cl(U)$ 
based at some point $\gamma(t_3)$ which is homotopic along 
$\gamma([t_2,t_3])$ to $C_1$. If $C_1$ is contractible in $M$, then the same
statement is true with $C_3$ equal to a constant curve.

Look at the universal cover $\tilde{M}$.  Let $\tilde{U}=\pi^{-1}(U)$
and $\tilde{V}=\pi^{-1}(V)$.  Let $\tilde{\gamma}$ be a lift of the ray 
$\gamma$ and $g\in \pi(M)$ be the deck transform represented by $[C_1]$.  It
may be the identity.  Let $\tilde{C_1}\in \tilde{U}$ be the lift of 
$C_1$ running from $\tilde{\gamma(t_2)}$ to $g \tilde{\gamma(t_2)}$ and
 $\tilde{C_3}\in \tilde{M}\setminus \tilde{U}$ be the lift of 
$C_3$ running from $\tilde{\gamma(t_3)}$ to $g \tilde{\gamma(t_3)}$. 

Then there exists $H:[0,1]\times [t_2,t_3] \longmapsto \tilde{M}$ 
such that
$H(s,0)=\tilde{C_2}(s)$, $H(s,1)=\tilde{C_3}(s)$, and 
$H(0,t)=\tilde{\gamma}(t)$ and $H(1,t)=g\tilde{\gamma(t)}$.  Here
we may have to reparametrize $C_2$ and $C_3$.  

Note that $H^{-1}(\tilde{U})$ and $H^{-1}(\tilde{V})$ are relatively 
open in $[0,1]\times [t_2,t_3]$ and their
union is $[0,1]\times [t_2,t_3]$.  
We would like to find a curve 
\begin{equation} \label{conncurve}
(s(r), t(r)) \subset 
H^{-1}(\tilde{U}) \cap H^{-1}(\tilde{V})\subset [0,1]\times [t_2,t_3]
\end{equation}
such that $h(0)=(0,t_2)$ and $h(1)=(1,t_2)$.  Then  
$$
C_2(r):=\pi(H(s(r), t(r))) \subset U \cap V
$$
is homotopic to $C_1$ based at $\gamma(t_2)$ and we are done.

To prove this we need only find a connected relatively open set contained in
$H^{-1}(\tilde{U}) \cap H^{-1}(\tilde{V})$ which contains both
$(0,t_2)$ and $(1,t_2)$.   This is true because connected open sets in
Euclidean space are pathwise connected. %See attachment to Feb 4

We employ the following lemma from Munkrees textbook [Mnk].

\noindent{\bf Lemma 13.1 of Munkrees}:
{\em
Let  $W=X\cup Y$ where $X$ and $Y$ are open sets and let
$X \cap Y =A \cup B$ where $A$ and $B$ are disjoint open sets.
If there exist two paths connecting $a \in A$ to $b\in B$, one
contained in $X$ and the other contained in $Y$ then $\pi_1(W,a)\neq 0$.}

We let $W=[0,1]\times [t_2,t_3]$, $X=H^{-1}(\tilde{U})$ and 
$Y=H^{-1}(\tilde{V})$.  Let $A$ be the connected component of
$X \cap Y$ which contains $a=(0,t_2)$ and let $B$ be the connected
component of $X \cap Y$ which contains $b=(1,t_2)$.  I claim $A=B$.
If not, they are disjoint and we can apply the lemma.  There exists a 
curve, namely $(s,t_2)$ contained in $X$ joining $a$ to $b$.  By our choice of
$t_2$ in (\ref{t2}), there exists a curve running around the other
three sides of the square joining $a$ to $b$ which is contained in $Y$.
So $\pi_1(W,a)\neq 0$.  This contradicts the fact that $W$ is contractible.
Thus $A=B$ is a connected component of
$H^{-1}(U) \cap H^{-1}(V)$ containing both $(0,t_2)$ and $(1,t_2)$.

\ProofEnd

\begin{theorem} \label{doublelocal}
Let $M^n$ be a complete Riemannian manifold with nonnegative Ricci
curvature.  Let $D\subset M$ be a precompact region
with smooth boundary, $\gamma$ a ray starting at $\gamma(0)\in D$ and 
$S$ be any  connected component of $\partial D$ containing a point $\gamma(t_1)$.
Then the image of the inclusion map
$$
i_*:\pi(S, \gamma(t_1)) \longmapsto \pi(Cl(D), \gamma(t_1))
$$
is $N \subset \pi(Cl(D), \gamma(t_1))$ such that $\pi(Cl(D), \gamma(t_1))/N $ 
contains at most two elements.

In fact, it contains only one element unless u 
\end{theorem}

\noindent{\bf Proof:}
If $M^n$ has the loops to infinity property, then by Theorem~\ref{Boundary},
we know that $\pi(Cl(D), \gamma(t_1))/N =\{e\}.$  If it does not, then
by Theorem~\ref{doublecover}, there is a double cover, $\bar{M}$, of $M^m$,
which splits along the lift, $\tilde{\gamma}$ of the geodesic, $\gamma$, 
and has the loops to infinty property.  

If $\pi^{-1}(Cl(D))$
is not connected, then $Cl(D)$ is homeomorphic to one of the connected
components of its lift.  So we can apply Theorem~\ref{Boundary} to the 
connected component, and we see that $i_*$ is a surjection.

If  $\pi^{-1}(Cl(D))$ is connected, then it is the double cover of
$Cl(D)$.  So there exists an element $g \in \pi_1(Cl(D), \gamma(t_1))$
whose representatives are lifted to nonclosed paths in $\pi^{-1}(Cl(D))$.  
We need only
show that if $h \in \pi_1(Cl(D), \gamma(t_1))$, then there exists 
$\bar{h}\in \pi(S, \gamma(t_1))$ such that either 
$i_*(\bar{h})=h$ or $i_*(\bar{h})=gh$.   

If $h \in \pi_1(Cl(D), \gamma(t_1))$ then either a representative lifts
to closed loops based at $\tilde{\gamma}(t_1))$
in $\pi^{-1}(Cl(D))$ or the representatives of $gh$ do.  Let $\tilde{C}$
be the lifted loop.  

By the loops to infinity property on the double cover and 
Theorem~\ref{Boundary}, there exists an element 
$\tilde{h}\in \pi_1(\pi^{-1}(S), \tilde{\gamma}(t_1))$
such that $[\tilde{C}]\in \tilde{h}$.  Let $\bar{h}=\pi_*(\tilde{h})$.
Then $\bar{h}\in \pi_1(S, {\gamma}(t_1))$ and 
$$
\bar{h}=[\pi \circ \tilde{C}]= h \textrm{ or } gh.
$$
\ProofEnd

\end{section}


\begin{thebibliography}{ChGrTa}


\bibitem[AbGl]{AbGl}  U. Abresch and D. Gromoll.     %Excess Top Est
 {\em On complete manifolds with nonnegative Ricci curvature.} 
 J. Amer. Math. Soc. 3 (1990), no. 2, 355--374.

\bibitem[And]{And} M. Anderson.
  {\em On the topology of complete manifolds of nonnegative Ricci
  curvature.} Topology 29 (1990), no. 1, 41--55.

\bibitem[doC]{doC} M.P. do Carmo. \underline{Riemannian Geometry}. 
 Translated by Francis Flaherty. Mathematics: Theory \& Applications. 
 Birkhauser Boston, Inc., Boston, MA, 1992. 

\bibitem[Ch]{Ch} J. Cheeger. {\em Critical points of distance functions and 
 applications to geometry.} Geometric topology: recent developments 
 (Montecatini Terme, 1990), 1--38, Lecture Notes in Math., 1504, Springer,
 Berlin, 1991. 

% soul thm
\bibitem[ChGl1] {ChGl1} J. Cheeger and D. Gromoll, {\em On the structure 
 of complete manifolds of nonnegative curvature.} Annals of Math. (2) 
 96 (1972), 413--443

\bibitem[ChGl2]{ChGl2} J. Cheeger and D. Gromoll  {\em
 The splitting theorem for manifolds of nonnegative Ricci curvature.} 
 J. Differential Geometry 6 (1971/72), 119--128. 

\bibitem[ChEb]{ChEb} J. Cheeger and D. Ebin. 
 \underline {Comparison Theorems in 
 Riemannian} \underline{ Geometry.} 
 North-Holland Mathematical Library, Vol. 9. 
 North-Holland Publishing Co., Amsterdam-Oxford; American
 Elsevier Publishing Co., Inc., New York, 1975. viii+174 pp. 

\bibitem[Fra] {Fra} T. Frankel,  {\em On the Fundamental Group of a
 Compact Minimal Submanifold.}  Ann. of Math. (2) 83 (1966) 68--73.

\bibitem[GlMy]{GlMy} D. Gromoll and W. Meyer, {\em  On complete open 
 manifolds of positive curvature.} Ann. of Math. (2) 90 (1969) 75--90. 

\bibitem[Law] {Law} H. B. Lawson, Jr, {\em The Unknottedness of Minimal
 Embeddings.}  Inventiones math. 11 (1970) 183-187.

\bibitem[Li] {Li} P. Li, \underline{Lecture Notes on Geometric Analysis}
 Lecture Notes Series, 6. Seoul National University, Research Institute 
 of Mathematics, Global Analysis Research Center, Seoul, 1993. 

%\bibitem[Li] {Li} P. Li, add ref

\bibitem[Mi]{Mi} J. Milnor, {\em A note on curvature and fundamental group.} 
 J. Differential Geometry 2 1968 1--7.

\bibitem[Mnk]{Mnk} James R. Munkres. \underline{Topology: a first course.}
  Prentice-Hall, Inc., Englewood Cliffs, N.J., 1975. xvi+413 %lemma 13.1

\bibitem[Na] {Na}  P. Nabonnand, {\em Sur les varietes riemanniennes 
 completes e courbure de Ricci positive.} (French) C. R. Acad. Sci. 
 Paris Sér. A-B 291 (1980), no. 10, A591--A593.

% ?
%\bibitem[Per] {Per} G. Perelman, {\em Collapsing with no proper extremal 
% subsets.} Comparison geometry (Berkeley, CA, 1993--94), 149--155, 
% Math. Sci. Res. Inst. Publ., 30, Cambridge Univ. Press, Cambridge, 1997. 

\bibitem[SchYau1] {SchYau1} R. Schoen and S-T Yau, {\em Complete 
 three-dimensional manifolds with positive Ricci curvature and scalar 
 curvature.} Seminar on Differential Geometry, pp. 209--228, 
 Ann. of Math. Stud., 102, Princeton Univ. Press, Princeton, N.J., 1982. 

\bibitem[SchYau2] {SchYau2} R. Schoen and S-T Yau, { Harmonic maps and 
 the topology of stable hypersurfaces and manifolds
 with non-negative Ricci curvature.}  Comment. Math. Helv. 51 (1976), 
 no. 3, 333--341.  %% CHECK REF OR USE HARM MAPS TEXT 

\bibitem [ShaYng] {ShaYng} J-P Sha and D-G Yang, {\em 
 Examples of manifolds of positive Ricci curvature.} 
 J. Differential Geom. 29 (1989), no. 1, 95--103. 

\bibitem[So] {So} C. Sormani, {\em Nonnegative Ricci Curvature, Small 
 Linear Diameter Growth, and Finite Generation of Fundamental Groups.}
 Preprint 1998.

\bibitem[Wei] {Wei}  G. Wei, {\em Examples of complete manifolds of positive 
 Ricci curvature with nilpotent isometry groups.}  Bull. Amer. Math. Soc. 
 (N.S.) 19 (1988), no. 1, 311--313.

\end{thebibliography}
\end{document}